\documentclass[11pt]{article}

\usepackage{amsmath,amsthm,amsfonts}

\usepackage{amssymb}

\usepackage{epsfig}

\usepackage{rotating}

\usepackage{subfigure}

\def\hR{{\hat{R}}}

\begin{document}

\title{Higher order Mori-Zwanzig models for the Euler equations}
\author{Panagiotis Stinis\\ 
\\
Department of Mathematics \\
    University of California \\
and \\
Lawrence Berkeley National Laboratory \\
    Berkeley, CA 94720}
\date {}

\maketitle

\begin{abstract}
In a recent paper \cite{CHSS06}, an infinitely long memory model (the t-model) for the Euler equations 
was presented and analyzed.  The model can be derived by keeping the zeroth order term in a Taylor expansion of the memory integrand in the Mori-Zwanzig formalism. We present here a collection of models for the Euler equations which are based also on the Mori-Zwanzig formalism. The models arise from a Taylor expansion of a different operator, the orthogonal dynamics evolution operator, which appears in the memory integrand. The zero, first and second order models are constructed and simulated numerically. The form of the nonlinearity in the Euler equations, the special properties of the projection operator used and the general properties of any projection operator can be exploited to facilitate the recursive calculation of even higher order models. We use our models to compute the rate of energy decay for the Taylor-Green vortex problem. The results are in good agreement with the theoretical estimates. The energy decay appears to be organized in "waves" of activity, i.e. 
alternating periods of fast and slow decay. Our results corroborate the assumption in \cite{CHSS06}, that the modeling of the 3D Euler equations by a few low wavenumber modes should include a long memory. 
\end{abstract}

\section{Introduction}
The advent of very powerful computers enhances our ability to probe complex systems and reveal 
their dynamics. Yet, there are many problems where the present computational capacity is not enough to fully resolve all their dynamic features. This situation poses a great challenge for numerical analysts. How can one extract meaningful information about a system's evolution when the direct simulation 
of such a system is out of reach? One approach is to use the qualitative information gathered by 
analytical, numerical and experimental studies as a guide for the formulation of reduced models 
which will hopefully reproduce the qualitative behavior observed so far, and moreover, reveal 
something new about the behavior of the system. Fluid turbulence is the archetypal example of 
a complex system and one of the main subjects of research for scientific computing. Even though a lot of information has been collected regarding the qualitative behavior of a real fluid (and of the associated equations of fluid flow) \cite{barenblatt,batchelor,bernard, C94, doering, kolmogorov,monin}, this information is very difficult to incorporate in the construction of reduced models. In the present paper we construct a collection of reduced models for the 3D Euler equations which describe the evolution of an inviscid fluid, based on the physical and numerical observation that 
the evolution of very smooth initial conditions can give rise to organized structures (known as vortices, vortex filaments, pancakes, sheets etc.) \cite{frisch}. These organized structures are known to create 
long temporal correlations \cite{alder} and since temporal correlations appear as memory integrands in our reduction formalism, the models we will construct are predominantly long memory ones.

The reduced models we will construct are based on the Mori-Zwanzig formalism of irreversible statistical mechanics \cite{mori,zwan1,zwan2} as reformulated by Chorin, Hald and Kupferman \cite{CHK00}. In this formalism, the correlations of the unresolved modes appear in the integrand of a memory term. In previous work \cite{CHK00}, it was found that the form of the memory term can simplify significantly if the approximation of a very long memory is made. This approximation can be effected by expanding in a Taylor series the memory integrand and keeping only the zeroth order term. That simplified model is known as the t-model and it was applied to the Euler equations  \cite{CHSS06} and to the Burgers equation  \cite{bernstein}. Here we propose an 
alternative construction which is also based on a Taylor expansion of an operator (the orthogonal dynamics operator) which appears in the memory term. The t-model can be derived as special case of the present construction. We construct reduced models that retain the zeroth, first and second order 
of the Taylor expansion of the orthogonal dynamics operator. In addition, we derive a set of rules that can facilitate the recursive calculation of high order models. The rules are based on the observation that the form of the terms appearing in the reduced models is determined by the form of the nonlinearity, the special properties of the projection operator used and the general properties of any projection operator.

It is known that the formation of organized structures in a fluid flow is manifested through a cascade of energy from the 
large to the small scales (see e.g. \cite{frisch}). Even though the Euler equations conserve the energy of a smooth solution, they do not have to conserve the energy of a non-smooth solution. There exist estimates (see \cite{speziale} and references therein) about the rate of decay of energy in an inviscid flow. We use our 
models to compute the rate of energy decay for the Taylor-Green vortex \cite{taylor,gottlieb}. The results are in good agreement with the theoretical estimates. It is interesting that the energy decay appears to be organized in "waves" of activity, i.e. alternating periods of fast and slow decay.

The paper is organized as follows. In Section \ref{eulerequations} we review the Euler equations and 
the problem of underresolved computations for these equations. In Section \ref{mori-zwanzig} we 
give a brief presentation of the Mori-Zwanzig formalism which is the starting point of our approximations. 
In Section \ref{taylor} we present the models based on the Taylor expansion of the memory term integrand. We also give a set of rules for the recursive evaluation of high order models. Numerical results of the models for the Taylor-Green vortex are presented in Section \ref{numbers}.  Finally, in Section \ref{conclusions} we discuss the numerical results and what they 
suggest for future work.

%%%%%%%%%%%%%End of Section 1%%%%%%%%%%%%%%%%%%

\section{The Euler equations and the problem of underresolved computations}\label{eulerequations}
Consider the 
3D incompressible Euler equations with periodic boundary conditions in the cube $[0,2\pi]^3$:

\begin{gather}
\label{euler}
v_t+v\cdot \nabla v= - \nabla p  ,\; \nabla \cdot v=0, 
\end{gather}
where $v(x,t)=(v_1(x_1,x_2,x_3,t),v_2(x_1,x_2,x_3,t),v_3(x_1,x_2,x_3,t))$ is the velocity, $p$ is the 
pressure and $\nabla= (\frac{\partial}{\partial x_1},\frac{\partial}{\partial x_2},\frac{\partial}{\partial x_3} ).$ The system in (\ref{euler}) is 
supplemented with the initial condition $v(x,0)=v_0(x)$ which is also periodic and incompressible and 
$x=(x_1,x_2,x_3).$

Since we are working with periodic boundary conditions, we expand the solution in Fourier series 
keeping $M$ modes in each spatial direction, 
$$v_{M}(x,t )=\underset{k \in F \cup G}{\sum} u_k(t) e^{ikx},$$
where $F \cup G=[-\frac{M}{2},\frac{M}{2}-1]\times[-\frac{M}{2},\frac{M}{2}-1]\times[-\frac{M}{2},\frac{M}{2}-1].$ Also $k=(k_1,k_2,k_3)$ and $u_k(t)=(u^1_k(t),u^2_k(t),u^3_k(t)).$ 
We have written the set of Fourier modes as the union of two sets 
in anticipation of the construction of the reduced model comprising only of the modes in $F=[-\frac{N}{2},\frac{N}{2}-1]\times[-\frac{N}{2},\frac{N}{2}-1]\times[-\frac{N}{2},\frac{N}{2}-1],$ where $ N < M.$
The equation of motion for the Fourier mode $u_k$ becomes
\begin{equation}
\label{eulerode}
 \frac{d u_k}{dt}=- i \underset{p, q \in F \cup G}{\underset{p+q=k  }{ \sum}} k \cdot u_{p} 
A_{k} u_{q},
\end{equation}
where $A_k= I - \frac{k k^T}{|k|^2}$ is the incompressibility projection matrix and $I$ is the $3\times3$ 
identity matrix. The symbol $\cdot$ denotes inner product in $\mathbb{R}^3.$ The system (\ref{eulerode}) is supplemented by the initial condition $u_0=\{u_k(0)\}=\{u_{0k}\}, \; k \in F \cup G,$ where $u_{0k}$ are the Fourier coefficients of the initial condition $v_0(x).$

Even if we start from a very smooth initial condition and $M$ is of the order of $10^3$ in each direction (the state of the art in massively parallel computers), the solution of the system of ordinary differential equations (\ref{eulerode}) can create significant activity in the highest modes of our allowed resolution. This phenomenon, called nonlinear instability, stems from the 
cascade (transfer) of energy from the large to the small scales and renders our calculations meaningless once a significant amount of energy has accumulated in the smallest scales of the solution. The source 
of the problem is that the system of equations (\ref{eulerode}) conserves the energy $E=\frac{1}{2}\underset{k \in F \cup G}{\sum} |u_k|^2.$ Thus, when energy reaches the highest modes in our calculation, there is 
no way for this energy to exit the range of allowed wavenumbers, so it begins to contaminate the 
solution and leads to a catastrophic increase of the error. Hence, if one hopes to extract meaningful 
information from a finite calculation, one needs to account for the required drain of energy out of 
the allowed range of wavenumbers (the same situation appears also in the case of the Navier-Stokes equations with small viscosity). Since we assume that we can only afford to solve the 
system (\ref{eulerode}) for $M$ modes in each direction (the set $F \cup G$), our task is to construct a reduced model for the modes in $F$ and use the modes in $G$ to effect the needed drain of 
energy out of the set $F.$ This is why we wrote the set of Fourier modes as a union of two sets. From 
now on, the set $F$ will be called the resolved modes and the set $G$ the unresolved modes.

A lot of work (e.g. \cite{foias,moser,pasquetti,piomelli,scotti,she}) has been devoted to 
the development of models that effect the drain of energy needed to keep the calculation 
well resolved. The problem with such models is that they are introduced {\it ad hoc}, involve parameters that need to be adjusted and usually work for some cases but not others. Another usual drawback is the perturbative treatment of the nonlinear term. While this is adequate for large values of the viscosity, it is not for the more realistic small viscosity cases (see \cite{smith} and references therein). There is an obvious need for 
models that are derived directly from the Euler equations based on assumptions that respect the 
observed physics. In order to do that one needs as a starting point a framework that allows, in principle, the construction of reduced models of arbitrary precision. After such a formalism has been established, 
reduced models of varying precision can be constructed by incorporating qualitative knowledge 
about the behavior of the solutions. One instance of such a formalism is the Mori-Zwanzig formalism 
of irreversible statistical mechanics \cite{zwan1,zwan2} as reformulated by Chorin, Hald and Kupferman 
\cite{CHK00,CH05}. As with every other formalism that allows the construction, in principle, of models 
of arbitrary accuracy, the actual formulation of the reduced model requires a lot of information about the 
full system. This could appear to defeat the purpose of constructing a reduced model. However, it 
is to be expected, because each system has a certain amount of information and a reduced model 
should account for this information. In \cite{CHSS06}, starting from the Mori-Zwanzig formalism and 
observations about the behavior of the solutions of fluids for very small or zero viscosity, an infinitely 
long memory model was proposed for the Euler equations. Such a model was based on the 
observation (see e.g. \cite{alder}) that the formation of vortices in a fluid flow can lead to 
very slowly decaying velocity temporal correlations. Since correlations are the building blocks of the 
memory term in the Mori-Zwanzig formalism, the assumption of an infinite long memory, which 
simplifies considerably the construction and the form of the reduced model, was natural. We should 
note here that the reduced model in \cite{CHSS06} (known as the t-model) goes against the 
common practice in modeling, where usually the opposite assumption  of extremely short, or none at all, 
memory is invoked. As was shown in \cite{CS05}, the two cases, of extremely long and 
extremely short memory, are the two sides of the same coin. However, there is a huge difference 
between them as far as their domain of validity. In the very short memory model, the unresolved modes 
are assumed to have very fast decaying temporal correlations, while the opposite is assumed in the case of the very long memory model. Such locality or non-locality of the correlations should ultimately be
determined by the choice of resolved modes \cite{shalizi03}.

%%%%%%%%%%%%%%End of Section 2%%%%%%%%%%%%%%%%%%%

\section{The Mori-Zwanzig formalism}\label{mori-zwanzig}

We give a brief account of the Mori-Zwanzig formalism (\cite{CHK00}) which is the starting 
point for the different models to be presented in the next section. Note that here we use 
the standard notation associated with the Mori-Zwanzig formalism, and the use of the 
variables $u,x$ should be clear from the context.

Suppose we are given an $M$-dimensional system of ordinary differential equations
\begin{equation}
\frac{d\phi}{dt}=R(\phi) \label{odes} 
\end{equation}
with initial condition $\phi(0)=x.$

The system of ordinary differential equations
we are asked to solve can be transformed into the linear
partial differential equation 
\begin{equation}
\label{pde}
u_t=Lu, \qquad u(x,0)=g(x)
\end{equation}
where $L=\sum_{i=1}^M R_i(x)\frac{\partial}{\partial{x_i}}$ and the solution of (\ref{pde}) is
given by $u(x,t)=g(\phi(x,t))$. Consider the following
initial condition for the PDE
$$g(x)=x_k \Rightarrow  u(x,t)=\phi_k(x,t)$$
Using semigroup notation we can rewrite (\ref{pde}) as
$$\frac{\partial}{\partial{t}} e^{tL}x_k=L e^{tL}x_k$$
Suppose that the vector of initial conditions can be divided as $x=(\hat{x},\tilde{x}),$ where 
$\hat{x}$ is the $N$-dimensional vector of the resolved variables and $\tilde{x}$ is the $(M-N)$-dimensional vector of the unresolved variables.  
Let $P$ be an orthogonal projection on the space of functions of $\hat{x}$ and $Q=I-P.$ In previous 
publications on the Mori-Zwanzig formalism, the projection operator $P$ was defined through 
a probability density function $f(x)$ on the set of variables $x.$ For the special case of the Euler equations, such a density is 
impossible to come by analytically or experimentally, thus we choose a different projection operator, 
which does not require the knowledge of a density for the values of $x.$ For a function $h(x)$ of all the 
variables, the projection operator we will use is defined by $P(h(\hat{x},\tilde{x}))=h(\hat{x},0),$ i.e. 
it replaces the value of the unresolved variables $\tilde{x}$ in any function $h(x)$ by zero. Similarly, the 
initial condition $x=(\hat{x},\tilde{x})$ is replaced by $(\hat{x},0).$ We should note here that such 
a projection operator is rather natural for the problem of the evolution of a very smooth initial condition, where one expects only a few Fourier modes to have nonzero values initially.  If we divide the 
wavenumbers in shells of different radii, then we can order the Fourier modes depending on which shell they belong to. For the numerical examples in this paper, the resolved modes will be taken in the 
first few shells, while the unresolved modes in the rest of the shells allowed by our resolution.

The equation (\ref{pde}) 
can be rewritten as 
\begin{equation}
\label{mz}
\frac{\partial}{\partial{t}} e^{tL}x_k=
e^{tL}PLx_k+e^{tQL}QLx_k+
\int_0^t e^{(t-s)L}PLe^{sQL}QLx_kds,
\end{equation}
where we have used Dyson's formula
\begin{equation}
\label{dyson1}
e^{tL}=e^{tQL}+\int_0^t e^{(t-s)L}PLe^{sQL}ds.
\end{equation}
Equation (\ref{mz}) is the Mori-Zwanzig identity. 
Note that
this relation is exact and is an alternative way
of writing the original PDE. It is the starting
point of our approximations. Of course, we
have one such equation for each of the resolved
variables $\phi_k, k=1,\ldots,N$. The first term in (\ref{mz}) is
usually called Markovian since it depends only on the values of the variables
at the current instant, the second is called "noise" and the third "memory". 
The meaning of the different terms appearing in (\ref{mz}) and a connection 
(and generalization) to the fluctuation-dissipation theorems of irreversible 
statistical mechanics can be found in \cite{CS05}. 

If we write
$$e^{tQL}QLx_k=w_k,$$ 
$w_k(x,t)$ satisfies the equation
\begin{equation}
\label{ortho}
\begin{cases}
&\frac{\partial}{\partial{t}}w_k(x,t)=QLw_k(x,t) \\ 
& w_k(x,0) = QLx_k=R_k(x)-PR_k(\hat{x}). 
\end{cases} 
\end{equation}
If we project (\ref{ortho}) we get
$$P\frac{\partial}{\partial{t}}w_k(x,t)=
PQLw_k(x,t)=0,$$
since $PQ=0$. Also for the initial condition
$$Pw_k(x,0)=PQLx_k=0$$
by the same argument. Thus, the solution
of (\ref{ortho}) is at all times orthogonal
to the range of $P.$ We call
(\ref{ortho}) the orthogonal dynamics equation. Since the solutions of 
the orthogonal dynamics equation remain orthogonal to the range of $P$, 
we can project the Mori-Zwanzig equation (\ref{mz}) and find
\begin{equation}
\label{mzp}
\frac{\partial}{\partial{t}} Pe^{tL}x_k=
Pe^{tL}PLx_k+
P\int_0^t e^{(t-s)L}PLe^{sQL}QLx_k ds.
\end{equation}

%%%%%%%%%%%%End of Section 3%%%%%%%%%%%%%%%%%%%%

\section{Mori-Zwanzig models based on the Taylor expansion of the orthogonal dynamics operator}\label{taylor}

The t-model \cite{CHK00,CHSS06} can be derived in different ways but the one that is closer to our approach consists of expanding the memory integrand $ e^{(t-s)L}PLe^{sQL}$ around $s=0$ and 
retaining only the zero order term. The essence of the t-model approximation is the approximation of the orthogonal dynamics operator $e^{tQL}$ by the full dynamics operator $e^{tL}.$ For the models presented here, we proceed in an 
alternative way by expanding the orthogonal dynamics operator around $s=0.$ Depending on how 
many terms we keep ($1,2,3,\ldots$), we obtain zeroth, first, second, .... order approximations respectively,
\begin{gather}
\label{orthoexp}
e^{sQL}=I+sQL+\frac{s^2}{2}QLQL+O(s^3), \\
PLe^{sQL}=PL+sPLQL+\frac{s^2}{2}PLQLQL+O(s^3).
\end{gather}
Every term in the expansion has one more factor of $QL$ than the previous term. At first sight, it seems 
that the construction of high order terms is just a tedious series of differentiations testing one's stamina. However, there are a few simple rules that can be used to facilitate the derivation of high order terms (see Section \ref{higher}).

There is no reason to expect {\it a priori} that the orthogonal dynamics operator can be 
written as a Taylor series of $s.$ In fact, we expect the dependence of $e^{sQL}$ on $s$ to be rather 
rough due to the rapid energy cascade to the unresolved modes.  The expression of this cascade is the fact that the orthogonal dynamics operator evolves quantities in a way that they acquire a component in the orthogonal complement of the range of the operator $P.$ However, the results obtained from models 
of different orders in $s$ are instructive and worth presenting (see also the discussion in Section \ref{conclusions}).

Before we present the different models, we rewrite the equations (\ref{eulerode}) to conform with 
the Mori-Zwanzig formalism. We set 
$$R_k(u)=- i \underset{p, q \in F \cup G}{\underset{p+q=k  }{ \sum}} k \cdot u_{p} A_{k} u_{q}$$
and we have 
\begin{equation}
\label{eulerodemz}
\frac{d u_k}{dt}=R_k(u) 
\end{equation}
for $ k \in F \cup G.$ The system (\ref{eulerodemz}) is supplemented by the initial 
condition $u_0=(\hat{u}_0,\tilde{u}_0)=(\hat{u}_0,0).$ Note that we focus on initial conditions where 
the unresolved Fourier modes are set to zero. This is enough, as mentioned before,  for our purposes, since the evolution of even very smooth initial conditions can still give rise to the phenomenon of 
nonlinear instability. We will proceed and construct reduced equations for the Fourier modes $u_k$ with $k \in F. $ Of course, these equations will depend on the values of the Fourier modes in $G$ and so, our task in constructing a model for the modes in $F,$ is to model the behavior of the modes in $G$ so that we can obtain a closed system for the modes in $F.$

\subsection{Zeroth, first and second order models}
We start our presentation of the models with the zeroth order model which is cubic in the Fourier modes.

\begin{equation}
\label{zeroa}
\frac{\partial}{\partial{t}} Pe^{tL}u_{0k}=Pe^{tL}\hR_k(\hat{u}_{0}) +P\int_0^t e^{(t-s)L} Z^0_k(\hat{u}_0) ds, \end{equation}
where
\begin{gather}
\label{zeroa2}
Z^0_k(\hat{u}_0)=
PLQLu_{0k}=-i   \biggl( \underset{p \in G, \;  q \in F}{\underset{p+q=k  }{ \sum}} k \cdot \hR_{p}(\hat{u}_0) 
A_{k} u_{0q} +  \underset{p \in F, \;  q \in G}{\underset{p+q=k  }{ \sum}} k \cdot u_{0p} 
A_{k} \hR_{q}(\hat{u}_0)  \biggr) 
\end{gather}
and 
$$\hR_k(\hat{u}_0)=R_k(\hat{u}_{0},0)=- i \underset{p, q \in F}{\underset{p+q=k  }{ \sum}} k \cdot u_{0p} A_{k} u_{0q}.$$
If we make the approximation $\int_0^t e^{(t-s)L} Z^0_k(\hat{u}_0) ds \sim  \int_0^t  e^{tL} Z^0_k(\hat{u}_0) ds=t e^{tL} Z^0_k(\hat{u}_0), $ then we recover the t-model of \cite{CHSS06}.

The system in (\ref{zeroa}) is not closed in the quantities $Pe^{tL}u_{0k}$ for $k \in F.$ The reason is that the projection operator $P$ does not in general commute with the {\it nonlinear} functions $e^{tL}R(\hat{u},0)$ and $e^{tL}Z^0(\hat{u},0).$ The simplest way to obtain a closed system is to commute the 
nonlinear functions and the projection $P.$ This is a kind of mean-field approximation. The difference 
with other mean-field approximations is that here we do account for the fluctuations of the unresolved modes, through the inclusion of the memory term. Thus, the reduced system becomes
\begin{equation}
\label{zerob}
\frac{\partial}{\partial{t}} Pe^{tL}u_{0k}=\hR_k(Pe^{tL}\hat{u}_{0}) +\int_0^t  Z^0_k(Pe^{(t-s)L}\hat{u}_0) ds, \end{equation}

The first order model is quartic in the Fourier modes and is given by
\begin{gather}
\label{firsta}
\frac{\partial}{\partial{t}} Pe^{tL}u_{0k}=\hR_k(Pe^{tL}\hat{u}_0)+\int_0^t  Z^0_k(Pe^{(t-s)L}\hat{u}_0)ds \\
+\int_0^t  s Z^1_k(Pe^{(t-s)L}\hat{u}_0) ds, 
\end{gather}
where
\begin{gather}
\label{firstb}
Z^1_k(\hat{u}_0)=
PLQLQLu_{0k}=  \notag \\
-i   \biggl( 
 \underset{p \in  F \cup G, \;  q \in G}{\underset{p+q=k  }{ \sum}} k \cdot \hR_{p}(\hat{u}_0) 
A_{k} \hR_{q}(\hat{u}_0)  +  \underset{p \in G, \;  q \in F \cup G}{\underset{p+q=k  }{ \sum}} k \cdot \hR_{p} 
(\hat{u}_0) A_{k} \hR_{q}(\hat{u}_0) + \notag \\
\underset{p \in G, \;  q \in F}{\underset{p+q=k  }{ \sum}} k \cdot Z^0_{p}(\hat{u}_0) 
A_{k} u_{0q} +  \underset{p \in F, \;  q \in G}{\underset{p+q=k  }{ \sum}} k \cdot u_{0p} 
A_{k} Z^0_{q}(\hat{u}_0)  \biggr)  \notag
\end{gather}

The second order model is quinitic in the Fourier modes and is given by
\begin{gather}
\label{seconda}
\frac{\partial}{\partial{t}} Pe^{tL}u_{0k}=\hR_k(Pe^{tL}\hat{u}_0)+\int_0^t  Z^0_k(Pe^{(t-s)L}\hat{u}_0)ds \\
+\int_0^t  s Z^1_k(Pe^{(t-s)L}\hat{u}_0) ds +\int_0^t  \frac{s^2}{2}Z^2_k(Pe^{(t-s)L}\hat{u}_0) ds, 
\end{gather}
where
\begin{gather}
\label{secondb}
Z^2_k(\hat{u}_0)=
PLQLQLQLu_{0k}= \notag \\
-i   \biggl( 
 \underset{p \in  F \cup G, \;  q \in G}{\underset{p+q=k  }{ \sum}} k \cdot Z^0_{p}(\hat{u}_0) 
A_{k} \hR_{q}(\hat{u}_0)  +  \underset{p \in G, \;  q \in F \cup G}{\underset{p+q=k  }{ \sum}} k \cdot \hR_{p} 
(\hat{u}_0) A_{k} Z^0_{q}(\hat{u}_0) + \label{secondb1} \\
 \underset{p \in  F \cup G, \;  q \in G}{\underset{p+q=k  }{ \sum}} k \cdot B_{p}(\hat{u}_0) 
A_{k} \hR_{q}(\hat{u}_0)  +  \underset{p \in G, \;  q \in F \cup G}{\underset{p+q=k  }{ \sum}} k \cdot \hR_{p} 
(\hat{u}_0) A_{k} B_{q}(\hat{u}_0) + \label{secondb2} \\
 \underset{p \in   G, \;  q \in F}{\underset{p+q=k  }{ \sum}} k \cdot Z^0_{p}(\hat{u}_0) 
A_{k} \hR_{q} (\hat{u}_0) +  \underset{p \in F, \;  q \in  G}{\underset{p+q=k  }{ \sum}} k \cdot \hR_{p}(\hat{u}_0)  A_{k} Z^0_{q}(\hat{u}_0) + \label{secondb3}\\
 \underset{p \in  F \cup G, \;  q \in G}{\underset{p+q=k  }{ \sum}} k \cdot Z^0_{p}(\hat{u}_0) 
A_{k} \hR_{q}(\hat{u}_0)  +  \underset{p \in G, \;  q \in F \cup G}{\underset{p+q=k  }{ \sum}} k \cdot \hR_{p} 
(\hat{u}_0) A_{k} Z^0_{q}(\hat{u}_0) + \label{secondb4}\\
 \underset{p \in  G, \;  q \in F \cup G}{\underset{p+q=k  }{ \sum}} k \cdot Z^0_{p}(\hat{u}_0) 
A_{k} \hR_{q}(\hat{u}_0)  +  \underset{p \in F \cup G, \;  q \in G}{\underset{p+q=k  }{ \sum}} k \cdot \hR_{p} 
(\hat{u}_0) A_{k} Z^0_{q}(\hat{u}_0) +\label{secondb5}\\
\underset{p \in G, \;  q \in F}{\underset{p+q=k  }{ \sum}} k \cdot Z^1_{p}(\hat{u}_0) 
A_{k} u_{0q} +  \underset{p \in F, \;  q \in G}{\underset{p+q=k  }{ \sum}} k \cdot u_{0p} 
A_{k} Z^1_{q}(\hat{u}_0)  \; \biggr) \label{secondb6}
\end{gather}
where 
$$B_k(\hat{u}_0)=- i \underset{p, q \in F}{\underset{p+q=k  }{ \sum}} k \cdot \hR_{p} (\hat{u}_0) A_{k} u_{0q}.$$

Some terms in (\ref{secondb1})-(\ref{secondb6}) can be grouped together to yield a simpler 
expression. Yet, it is more instructive to keep them separate because they reveal the 
rules that determine the types of expressions that appear in a term of arbitrary order.

One common aspect of all the models of order one and higher, is that they involve integrodifferential equations where 
the integrals are convolutions. The computation of convolution integrals can be very expensive. However, the form of the convolution integral in our case can be translated into a sum of ordinary integrals by a simple change of variables. For the first order model (\ref{firsta}) the change of 
variables $s'=t-s$ gives
\begin{gather*}
\frac{\partial}{\partial{t}} Pe^{tL}u_{0k}=\hR_k(Pe^{tL}\hat{u}_0)+\int_0^t  Z^0_k(Pe^{sL}\hat{u}_0)ds+\int_0^t  (t-s) Z^1_k(Pe^{sL}\hat{u}_0) ds= \\
\hR_k(Pe^{tL}\hat{u}_0)+\int_0^t  Z^0_k(Pe^{sL}\hat{u}_0)ds +t \int_0^t   Z^1_k(Pe^{sL}\hat{u}_0) ds -\int_0^t  s Z^1_k(Pe^{sL}\hat{u}_0) ds ,
\end{gather*}
The resulting integrals are no longer of convolution type and need not be evaluated from scratch for each value of $t.$ Instead, they can be evaluated by adding the 
contribution of each new timestep to the existing value of the integral. Similar constructions can be 
used for the convolution integrals appearing in the higher order terms.

Instead of expanding the orthogonal dynamics operator, one can start with equation (\ref{mzp}) 
\begin{equation*}
\frac{\partial}{\partial{t}} Pe^{tL}u_{0k}=
Pe^{tL}PLu_{0k}+
P\int_0^t e^{sL}PLe^{(t-s)QL}QLu_{0k} ds.
\end{equation*}
(note the change of variables $s'=t-s$) and build an infinite hierarchy of {\it Markovian} equations which is equivalent to the non-Markovian 
equation (\ref{mzp}). Truncating the hierarchy at some term is equivalent to keeping terms up to 
a certain order in the expansion of $e^{tQL}.$ In particular, if we define $w_{k0}(\hat{u}_0,t)=P\int_0^t e^{(t-s)L}PLe^{sQL}QLu_{0k} ds,$ (where we have denoted explicitly the dependence of the quantity $w_{k0}$ on the resolved initial conditions  $\hat{u}_0$), we get
$$\frac{dw_{k0}}{dt}=Pe^{tL}PLQLu_{0k}+P\int_0^t e^{sL}PLe^{(t-s)QL}QLQLu_{0k}ds.$$
We can define $$w_{k1}(\hat{u}_0,t)=P\int_0^t e^{sL}PLe^{(t-s)QL}QLQLu_{0k}ds,$$ then 
$$\frac{dw_{k1}}{dt}=Pe^{tL}PLQLQLu_{0k}+P\int_0^t e^{sL}PLe^{(t-s)QL}QLQLQLu_{0k}ds.$$
Similarly, we can define terms $w_{k2},\ldots,w_{kn},\ldots.$ For the $n$th quantity $w_{kn}$ we have
 $$\frac{dw_{kn}}{dt}=Pe^{tL}PL(QL)^nQLu_{0k}+P\int_0^t e^{sL}PLe^{(t-s)QL}(QL)^{n+1}QLu_{0k}ds.$$ Of course, we have to truncate the hierarchy at some point. This amounts to setting the integral term appearing in the equation of evolution for $w_{kn_0}$ (for some index $n_0$) equal to zero. Gathering the equations up to term $n_0$ we get
\begin{align*}
\frac{d}{dt}Pe^{tL}u_{0k}&=Pe^{tL}PLu_{0k}+w_{k0}(t),\\
\frac{dw_{k0}}{dt}&=Pe^{tL}PLQLu_{0k}+w_{k1}(t), \\
\frac{dw_{k1}}{dt}&=Pe^{tL}PLQLQLu_{0k}+ w_{k2}(t), \\
\ldots \\
\frac{dw_{kn_0}}{dt}&=Pe^{tL}PL(QL)^{n_0}QLu_{0k}
\end{align*}
By commuting the projection $P$ with the nonlinear functions appearing on the RHS of the equations, 
the above system of equations can be transformed into a closed system for $Pe^{tL}u_{0k}.$

\subsection{Higher order models}\label{higher}

Higher than the second order terms in the expansion of the orthogonal dynamics operator are hard 
to derive without some kind of book keeping rules, which can help reduce the risk of a mistake. There are three determining factors in the appearance of the expressions for a term of some order. 

\begin{enumerate}
\begin{item}
The form of the nonlinearity,
\end{item}

\begin{item}
The specific form of the projection,
\end{item}

\begin{item}
The general property of any projection $PQ=P(I-P)=0.$
\end{item}
\end{enumerate}
We proceed to analyze how the three factors mentioned above determine the form of the expressions appearing in higher order terms. This will lead us eventually to the formulation of a general scheme 
for computing higher order terms. First, the quadratic nonlinearity in the Euler equations and the fact that each new term in the series invloves one more differentiation (the $L$ in the factor $QL$), result in 
the term of order $n$ involving $n+3$ powers of the Fourier modes. For the zeroth, first and second order models this gave cubic, quartic and quintic powers as already pointed out. Also, the $n$th order 
term will involve all combinations of $n+3$ of the form $n+3=m+l,$ where $m,l$ are positive integers. For example, the second order term ($n=2$) involves quintic expressions ($n+3=5$) of the form $4+1,1+4,3+2,2+3$. In addition, the cubic expressions can be further decomposed as $2+1,1+2,$ (see $B_k(\hat{u}_0)$) and, of course, the trivial decomposition $3=3+0=0+3$ which is nothing but $Z^0(\hat{u}_0).$ On the other hand, the quartic expression does not appear in a decomposed form but only as $Z^1(\hat{u}_0).$ So, there is something more to 
the way the higher order terms than just a combination of all the possible arithmetic combinations of $n+3.$ In fact, there is much more structure that comes from the next two factors above, namely the special properties of the projection used and also the general property of every projection operator. However, before we go into more details let us make a comment: the structure of the expressions is convenient for numerical implementation because the different expressions are convolution sums in Fourier space, which can be computed through the Fast Fourier Transform (FFT) in real space, by FFT of appropriate arrays. The adjective appropriate is the catch. Even though all the expressions are convolution sums, one has to find what are the ranges for the wavenumber indices appearing in the sums. Factors 2) and 3) will help us determine what is this range.

To reveal this structure, we should examine more closely how a term in the series is related to the 
one preceding it. The $n$th order term has the form $Z^n(\hat{u}_0)=PL(QL)^nQLu_{0k},$ i.e. $n$ applications of 
the operator $QL$ and then application of the operator $PL.$ We have $PL(QL)^nQLu_{0k}=PLQL(QL)^{n-1}QLu_{0k}.$ The part $(QL)^{n-1}QLu_{0k}$ is common with the $n-1$st term  $Z^{n-1}(\hat{u}_0)=PL(QL)^{n-1}QLu_{0k}.$ When we act on $(QL)^{n-1}QLu_{0k}$ with the extra factor $QL=L-PL,$ we get
$QL(QL)^{n-1} QLu_{0k}=L(QL)^{n-1} QLu_{0k}-Z^{n-1}(\hat{u}_0).$ This is an important point. It tells 
us that the expression for $QL(QL)^{n-1} QLu_{0k}$ contains 3 types of terms: i) Terms that could not 
appear in $Z^{n-1}$ because of the special property of the projection which sets to zero expressions 
linear in $u_{0k}$ for $k \in G;$ ii) Terms that could not appear in $Z^{n-1}$ due to the general 
property of any projection that $PQ=0$ and iii) Terms of the form $h(u_0)-(Ph)(\hat{u}_0),$ where $(Ph)(\hat{u}_0)$ is any expression appearing in the term $Z^{n-1}.$ So, we can assemble the expressions 
appearing in $QL(QL)^{n-1} QLu_{0k}$ into three groups according to i),ii) and iii). Then we can apply the operator $PL$ once and we are done. Why is this grouping of expressions helpful?

Let us start with the expressions iii). For those expressions, the final application of the operator $PL$ needed to complete the calculation of $Z^n$ is trivial. For example, consider the case of the expression 
$$\underset{p \in G, \;  q \in F}{\underset{p+q=k  }{ \sum}} k \cdot Z^0_{p}(\hat{u}_0) A_{k} u_{0q}=\underset{p \in G, \;  q \in F}{\underset{p+q=k  }{ \sum}} k \cdot PLQLu_{0p}(\hat{u}_0) A_{k} u_{0q}$$ which is of type $3+1$ and appears in the term $Z^1$ (we have omitted the $-i$ factor).  This will give rise to a type iii) expression of the form 
\begin{gather*}
\underset{p \in G, \;  q \in F}{\underset{p+q=k  }{ \sum}} k \cdot LQLu_{0p}({u}_0) A_{k} u_{0q}-\underset{p \in G, \;  q \in F}{\underset{p+q=k  }{ \sum}} k \cdot PLQLu_{0p}(\hat{u}_0) A_{k} u_{0q}= \\
\underset{p \in G, \;  q \in F}{\underset{p+q=k  }{ \sum}} k \cdot QLQLu_{0p}({u}_0) A_{k} u_{0q}.
\end{gather*}
Application of $PL$ on this is trivial, since only the expression resulting  from the application of $L$ on  $QLQLu_{0p}({u}_0)$ can survive the action of $P.$ This is because the expression which comes from acting with $L$ on $ u_{0q},$ i.e. $$\underset{p \in G, \;  q \in F}{\underset{p+q=k  }{ \sum}} k \cdot QLQLu_{0p}({u}_0) A_{k} R_{q}(u_0)$$ will be killed by the action of $P,$ since $PQ=0.$ But then what 
remains from the action of $PL$ is
$$\underset{p \in G, \;  q \in F}{\underset{p+q=k  }{ \sum}} k \cdot PLQLQLu_{0p}({u}_0) A_{k} u_{0q}=
\underset{p \in G, \;  q \in F}{\underset{p+q=k  }{ \sum}} k \cdot Z^1_p({u}_0) A_{k} u_{0q},$$
a $4+1$ term.

For the types i) and ii) there is a little more work to do. We begin with 
$$ \underset{p \in  F \cup G, \;  q \in G}{\underset{p+q=k  }{ \sum}} k \cdot LR_{p}({u}_0) 
A_{k} u_{0q} $$
which is a $3+1$ type i) term that appears in $QLQLQLu_{0k}$ during the construction of $Z^2.$ This 
term could not have appeared in the term $Z^1$ due to the special property of our projection $P$ 
which sets to zero $u_{0q}$ with $q \in G.$ This term arises from the application of $L$ on a $2+1$ term. 
It is obvious that when acted upon with $PL$ it will give rise to a $3+2$ term only, since the associated $4+1$ term would vanish, again by the special property of the projection. Thus, it will contribute the 
terms 
\begin{gather*}
 \underset{p \in  F \cup G, \;  q \in G}{\underset{p+q=k  }{ \sum}} k \cdot (Z^0_{p}(\hat{u}_0) +B_{p}(\hat{u}_0) )A_{k} \hR_{q}(\hat{u}_0)  =\\
 \underset{p \in  F \cup G, \;  q \in G}{\underset{p+q=k  }{ \sum}} k \cdot Z^0_{p}(\hat{u}_0) 
A_{k} \hR_{q}(\hat{u}_0)  + 
 \underset{p \in  F \cup G, \;  q \in G}{\underset{p+q=k  }{ \sum}} k \cdot B_{p}(\hat{u}_0) 
A_{k} \hR_{q}(\hat{u}_0) 
\end{gather*}
where we have used the fact that $(PLR_{p})(\hat{u}_0)= Z^0_{p}(\hat{u}_0) +B_{p}(\hat{u}_0).$ This 
is nothing more but a result of the possible decompositions of the number $3$ as $3+0$ or $2+1,$ i.e. 
$ Z^0_{p}(\hat{u}_0)$ and  $B_{p}(\hat{u}_0).$ So, we see that the $3+1$ term which came from 
a $2+1$ term and the special property of the projection, gives rise to a $3+2$ term.

Let us examine a type ii) term. Consider the term 
$$ \underset{p \in G, \;  q \in F }{\underset{p+q=k  }{ \sum}} k \cdot (R_p- \hR_{p}) 
({u}_0) A_{k} R_{q}({u}_0) =  \underset{p \in G, \;  q \in F }{\underset{p+q=k  }{ \sum}} k \cdot QLu_{0p} A_{k} R_{q}({u}_0)$$
which is a $2+2$ type ii) term that appears in $QLQLQLu_{0k}$ during the construction of $Z^2.$ This 
term could not have appeared in $Z^1$ because of the general property of any projection $P$ that $PQ=0.$ It came from the application of $L$ on a $2+1$ term. It is obvious that when acted upon with 
$PL,$ the $2+2$ term will give rise to a $3+2$ term, since the associated $2+3$ term would vanish, again by the property $PQ=0.$ Thus, it will contribute the term
$$\underset{p \in G, \;  q \in F }{\underset{p+q=k  }{ \sum}} k \cdot PLQLu_{0p} A_{k} \hR_{q}(\hat{u}_0)=
\underset{p \in G, \;  q \in F }{\underset{p+q=k  }{ \sum}} k \cdot Z^0_p(\hat{u}_0) A_{k} \hR_{q}(\hat{u}_0)$$

In summary, we can write down the following rules for the evaluation of the $n$th order term $Z^n(\hat{u}_0)=PL(QL)^nQLu_{0k}$ in the Taylor series of the orthogonal dynamics operator:

\begin{enumerate}

\begin{item}
Write down the expression $(QL)^{n-1}QLu_{0k}.$
\end{item}
\begin{item}
Apply $QL$ and assemble the terms in the expression $QL(QL)^{n-1}QLu_{0k}$ in 3 groups: i) Terms that could not appear in $Z^{n-1}$ because of the special property of the projection which sets to zero expressions linear in $u_{0k}$ for $k \in G;$ ii) Terms that could not appear in $Z^{n-1}$ due to the general property of any projection that $PQ=0$ and iii) Terms of the form $h(u_0)-(Ph)(\hat{u}_0),$ where $(Ph)(\hat{u}_0)$ is any expression appearing in the term $Z^{n-1}.$
\end{item}

\begin{item}
Apply the operator $PL$ to the type i) terms. An $m+1$ term in  the expression $QL(QL)^{n-1}QLu_{0k}$ arose 
from an $(m-1)+1$ term and will give rise to an $m+2$ term. The symmetric term $2+m$ should 
also appear (the symmetric terms appear due to the rule of differentiating a product.)
\end{item}

\begin{item}
Apply the operator $PL$ to the type ii) terms. An $m+l$ term in  the expression $QL(QL)^{n-1}QLu_{0k}$ arose 
from an $m+(l-1)$ term where the $m$ term is of the form $(Qh)(u_0)$ for some function $h(u_0).$ It will give rise to an $(m+1)+l$ term. The symmetric $l+(m+1)$ term should also appear.
\end{item}

\begin{item}
Apply the operator $PL$ to the type iii) terms. There are two cases: a) An $(n-2)+1$ term will give rise to an $(n-1)+1$ term with the $n-1$st part being equal to $Z^{n-1}.$ The symmetric term $1+(n-1)$ should 
also appear; b) An $m+l$ term with $l \neq 1$ will give rise to an $(m+1)+l$ and an $m+(l+1)$ term. The symmetric terms $l+(m+1)$ and $(l+1)+m$ should also appear. 
\end{item}

\begin{item}
Make sure that in the final expressions {\it all} possible decompositions of $n+3$ into sums of two positive 
integers appear. All the expressions for the $n$th term in the series should be $n+3$ powers of 
Fourier modes.
\end{item}

\begin{item}
As a last resort, forget about the rules and proceed with straightforward differentiation. 
\end{item}

\end{enumerate}

We return to the expressions (\ref{secondb1})-(\ref{secondb6}) for $Z^2.$ Based on the rules 
above and the form of $Z^1,$ one can see that the expressions (\ref{secondb1})-(\ref{secondb2}) 
arose from type i) terms, the expression (\ref{secondb3}) arose from a type ii) term and the expressions (\ref{secondb4})-(\ref{secondb6}) arose from type iii) terms. It is evident that the proliferation of expressions is rather rapid and, even with the rules, the amount of work to derive high order terms can quickly become significant.

Note that the expansion of the orthogonal dynamics operator $e^{sQL}$ and,  more generally, of the expression $PLe^{sQL}$ in a Taylor series is equivalent to expanding the 
response function of the orthogonal dynamics to the "field" created by the resolved modes. Thus, 
the expressions for the terms in the Taylor series of the $PLe^{sQL}$ are what is known in the physics literature as sum rules \cite{forster}.

One final comment about the different models. By construction, all the reduced models are 
incompressible. This is because the terms appearing in the models involve the incompressibility projection operator $A_k.$

%%%%%%%%%%%%%%%%%%%%End of Section 4%%%%%%%%%%%%%%%%%

\section{Numerical results for the Taylor-Green problem}\label{numbers}

In this section we present numerical results of the application of the different models to the 3D Euler equations using the Taylor-Green vortex as an initial condition. 

\subsection{Zeroth order model}\label{zero}

Note that due to the quadratic nonlinearity, the equation for each Fourier mode in the Euler equations contains interactions with Fourier modes of at most double the wavevector. This means that if 
one models the Fourier modes with wavevectors at most double than the resolved modes, then 
one obtains a closed system of equations for the resolved modes. This means that the 
ratio of the number of modes included in the set $G$ over those in the set $F$ should be at least 1. Of course, 
one can have a larger set $G,$ depending on the computational power at hand. The more 
modes included in $G$, the better the modeling of the unresolved Fourier modes. The larger the ratio of the number of modes in $G$ over those in $F,$ the better the reduced model for the modes in $F$  should become.  We present results only for the case when the ratio is 1.

The different terms appearing in the RHS of the equations for the reduced models can be computed in 
real space using FFTs of appropriate arrays. Since for a reduced model of size $N$ in each spatial direction we include $N$ additional unresolved modes in each direction, the arrays involved in the FFTs should be of size $2N$. The fact that the model terms can be computed using the FFT makes their numerical 
implementation computationally efficient. Moreover, for the zeroth order model and for some expressions in the higher order terms of the form 
$$\underset{p \in G, \;  q \in F }{\underset{p+q=k  }{ \sum}} k \cdot H_p A_{k} C_{q},$$ the FFT calculations involved 
are dealiased by construction and thus no extra (e.g. $3/2$ rule \cite{canuto}) dealiasing is needed. This is straightforward to see. For a calculation involving $N$ modes in each direction, i.e. $N/2$ positive and $N/2$ negative, we perform FFTs of size $2N$, i.e. $N$ positive and $N$ negative modes. But we are 
interested only on the RHS for the first $N/2$ modes. This means (see \cite{canuto}) that to avoid 
aliasing (in the zeroth order term) we need for the total number of modes M used to satisfy the following inequality: $-N-N/2 \geq N/2-1-M$ which yields $M \geq 2N.$ But $2N$ is exactly how many modes 
we use in the FFTs, and thus the zeroth model term calculation through FFTs is dealiased by construction.

We have shown in Section \ref{taylor} that the memory term can be decomposed into a sum of 
ordinary integrals. This can make the calculation of the memory very efficient. Unfortunately, when 
implemented, the models need to have the range of integration for the memory term reduced from $[0,t]$ to $[t_0,t],$ otherwise the calculation becomes unstable. The value of $t_0$ is dependent on the 
model and the initial condition used. However, there is no tuning needed. The results become better, the longer the range of integration, until the value of the range that leads to instability. Thus, trial and error 
is needed not to fit the results to some prescribed curve, but just to find when does the calculation becomes unstable. A Taylor series around the current instant cannot be expected 
to be accurate for long times in the past and this is the reason for the need to truncate the memory term's range of integration. Even though we have to truncate the range of integration, we can still salvage some 
of the efficiency gained by the transformation of the convolution integral into a sum of ordinary 
integrals. Let us assume that we want to calculate the integral 
$$\int_{t_k+\Delta t- t_0}^{t_k+ \Delta t} f(s) ds,$$
where $f(s)$ is any of the integrands appearing in the different order models' memory terms. Decompose the integral as
$$\int_{t_k - t_0}^{t_k} f(s) ds +\int_{t_k}^{t_k+ \Delta t} f(s) ds -\int_{t_k- t_0}^{t_k+ \Delta t-t_0} f(s) ds$$
The first integral is already computed at the last step. The second integral is the contribution of the 
current step. And the third integral is the contribution that needs to be subtracted due to the 
fact that the memory has a truncated range of integration. Note that if the integral is not truncated, i.e. $t_0=t_k+\Delta t$ then the above decomposition is equal to 
$$\int_{0}^{t_k} f(s) ds +\int_{t_k}^{t_k+ \Delta t} f(s) ds$$
This shows that if there is no truncation, then the computation of the memory term is very efficient, since 
one needs only to add the contribution to the memory from the current step. 
We see that the truncation of the integral forces us to keep track of the values of the integrand at the instants 
$t_k+\Delta t -t_0$ and $t_k-t_0.$ But this means, that as the calculation progresses, we need to keep 
an array of length $[t_0/\Delta t],$ where $[]$ stands for integer part. This array needs to be updated 
at the end of every step so that it always keeps the values of the integrand for the last $[t_0/\Delta t]$ steps. At every step we will use only two values from the array, but we need to keep the whole 
history of length $[t_0/\Delta t]$ because, the range of integration extends only for $t_0$ to the past. Even though the integral does {\it not} need to be computed from scratch at every step, as 
for the case of a convolution integral, the process of updating the array of past integrand values can be expensive for 3D calculations. Also, the size of the array is considerable. For example, if we keep $N$ modes in each direction, the size of the array is $3N^3\times [t_0/\Delta t]$ (the factor of 3 comes from the 3 components of the velocity). Still, the computational time is half the one required if we keep the 
convolution integral form for the memory term.

In order to study the asymptotic decay rate of the energy in the resolved modes, one has to evolve the system for long times. We evolved each case up to time $t=100$, so that we have enough points to perform an accurate estimate of the decay rate exponent. The equations of motion for the Fourier modes were solved by the
modified Euler method  \cite{hair}. The integrals for the memory term were computed with the 
trapezoidal rule. The stepsize was set to $\Delta t=10^{-3}.$ We also performed an experiment with the Runge-Kutta 4th order method (again $\Delta t=10^{-3}$) and Simpson's rule for the evaluation of the integrals. The results were 
practically the same as with the lower order method (the difference between the results was not larger than $10^{-9}$). So, we decided to perform the rest of the experiments with the lower order method. 
As we mentioned before, the need to keep the values of the integrands in the past can be expensive. 
For our choice of stepsize and choice of integration interval for the memory term, we could only afford to use $8^3$ resolved modes when implementing the model on a single processor workstation.

It is important that the solutions respect the incompressibility condition. We know that the models are 
incompressible by construction. The calculation of the integral terms involves the summation of many 
instants of the solution. Even though for short times, the solution at each instant is incompressible to double precision, eventually, the inaccuracy in the summation involved in the integral terms can cause the solution to start deviating from incompressibility. However, for our experiments, even after $10^5$ steps the divergence never became larger than $10^{-14}.$ One can try to enforce the incompressibility condition by performing a projection on divergence-free fields at the end of each step. We tried that and the results did not change from the case without such a projection. This means that the violation of the divergence-free condition observed in our experiments is harmless.

\begin{figure}
\centering
%\subfigure[]{\epsfig{file=fig_zero_energy.eps, width=2.in}}
\subfigure[]{\epsfig{file=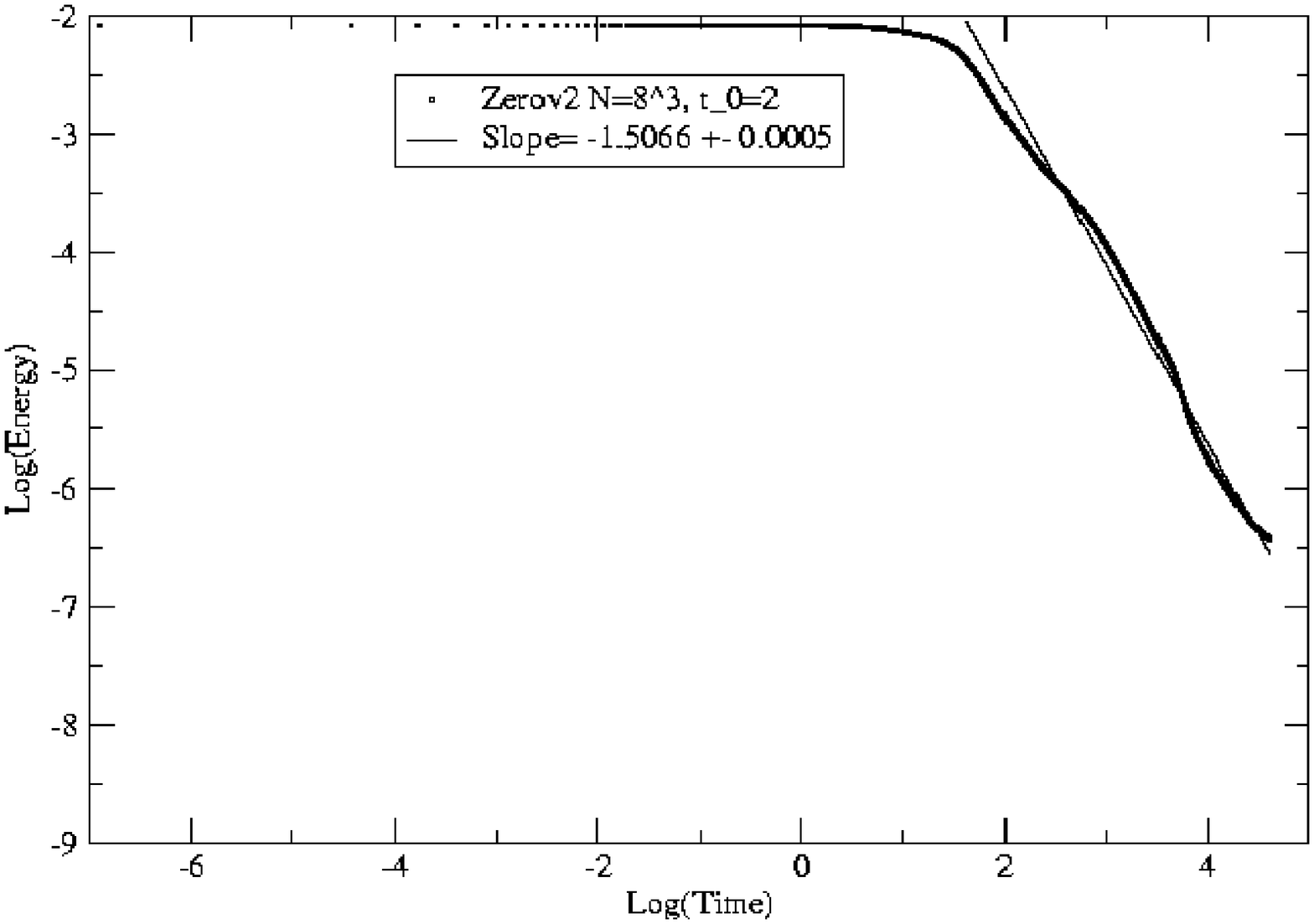, width=2.in}}
\qquad
%\subfigure[]{\epsfig{file=fig_zero_rate.eps,width=2.in}}
\subfigure[]{\epsfig{file=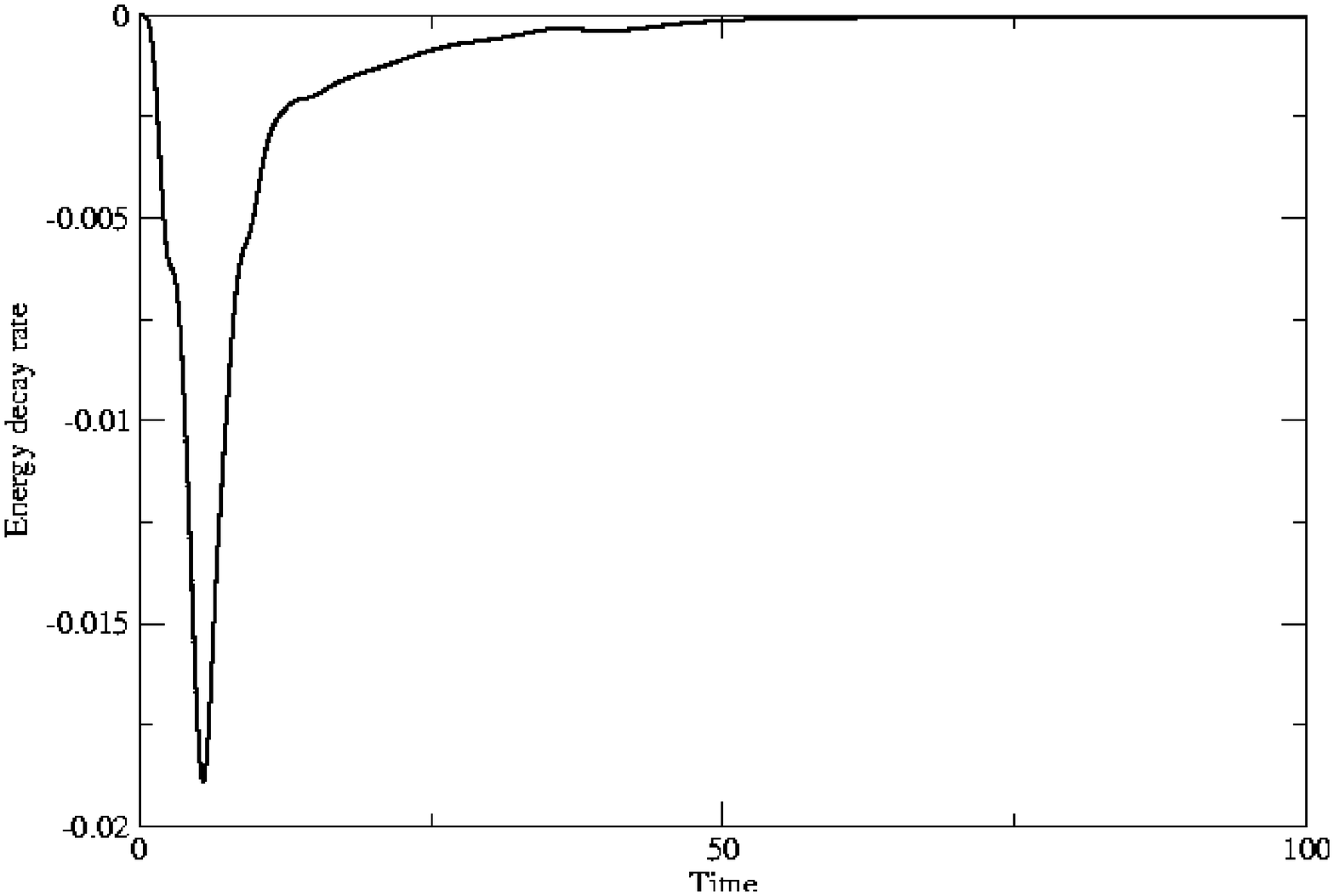,width=2.in}}
\caption{(a) Energy evolution for the zeroth order model with $N=8^3$ modes. (b) Evolution of the energy decay rate. }
\label{fig_zero}
\end{figure}

We present in Figure \ref{fig_zero}(a) the evolution of the energy in the resolved modes $E= \frac{1}{2}\underset{k \in F}{\sum} |u_k|^2$ for the zeroth order model with $8^3$  resolved modes. The 
slope (in log-log coordinates) is $\alpha=-1.5066 \pm 0.0005,$ which means that the energy decays as $t^{\alpha}.$ In Figure \ref{fig_zero}(b) we present results for the rate of energy decay $dE/dt$ for the zeroth order model with $8^3$ resolved modes. We set the parameter $t_0$ to the value 2. This is the longest history we are allowed to keep to avoid 
an instability.

The energy evolution shows an interesting trend. The energy decay seems to 
be organized in "waves" of activity, alternating periods of fast and slow decay. This behavior is absent in the t-model for the current resolution. Recall that the t-model is also a zeroth order model, albeit for the expansion of the whole memory integrand and not of the orthogonal dynamics operator only as the models here. In that paper, we found that one has to increase the number of resolved modes before the 
trend appears (more about the comparison between the two methods below).

This organization of the energy decay is reminiscent of the phenomenon of intermittency, i.e. bursts of activity followed by 
intervals of relative inaction on the part of the flow. Of course, the phenomenon of intermittency is not 
only of temporal nature, but has a spatial manifestation too. This is exhibited as concentration of the 
highest vorticity in small regions of the flow. The trend we observe in the decay of the energy seems to assign a specific purpose to the vorticity. Starting from a smooth initial condition, we have a steepening of the gradients in the field. This means that smaller scales are excited until the vorticity producing 
mechanism runs out of steam. Then we enter a period of relative inaction, until there is a restart of the 
mechanism of steepening. Energy is transferred again to the smaller scales and so forth. This scenario continues until there is no energy left in the large scales. After that, the flow just disintegrates and eventually comes to a halt. The purpose of vorticity mentioned above is to regulate the transfer of energy to the small scales. This is reminiscent of the picture suggested by Moffatt, Kida and Okhitani \cite{moffatt} of the vortex structures  acting as the "sinews of turbulence" .

\subsection{First order model}

\begin{figure}
\centering
%\subfigure[]{\epsfig{file=fig_first_energy.eps,width=2.in}}
\subfigure[]{\epsfig{file=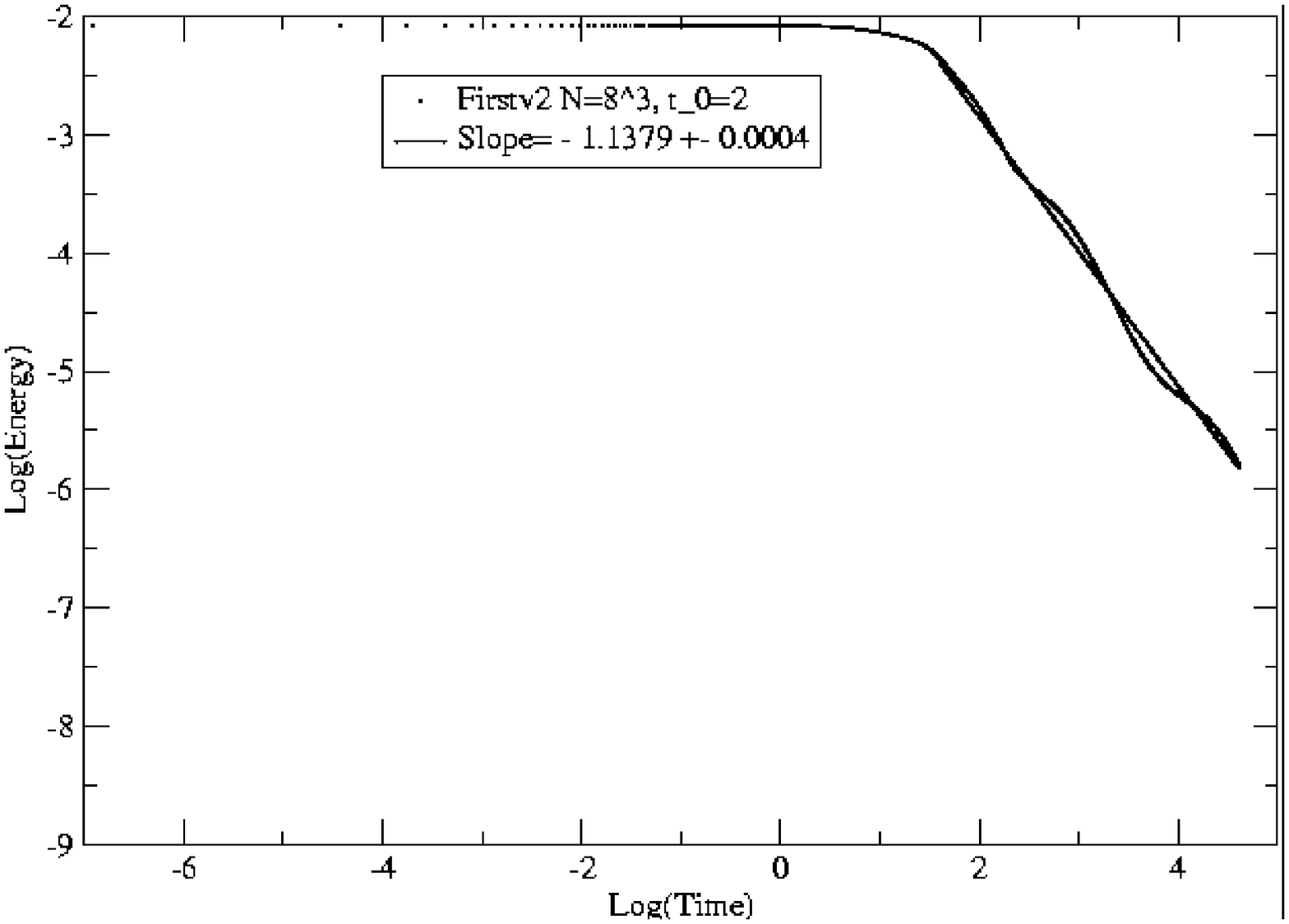,width=2.in}}
\qquad
%\subfigure[]{\epsfig{file=fig_first_rate.eps,width=2.in}}
\subfigure[]{\epsfig{file=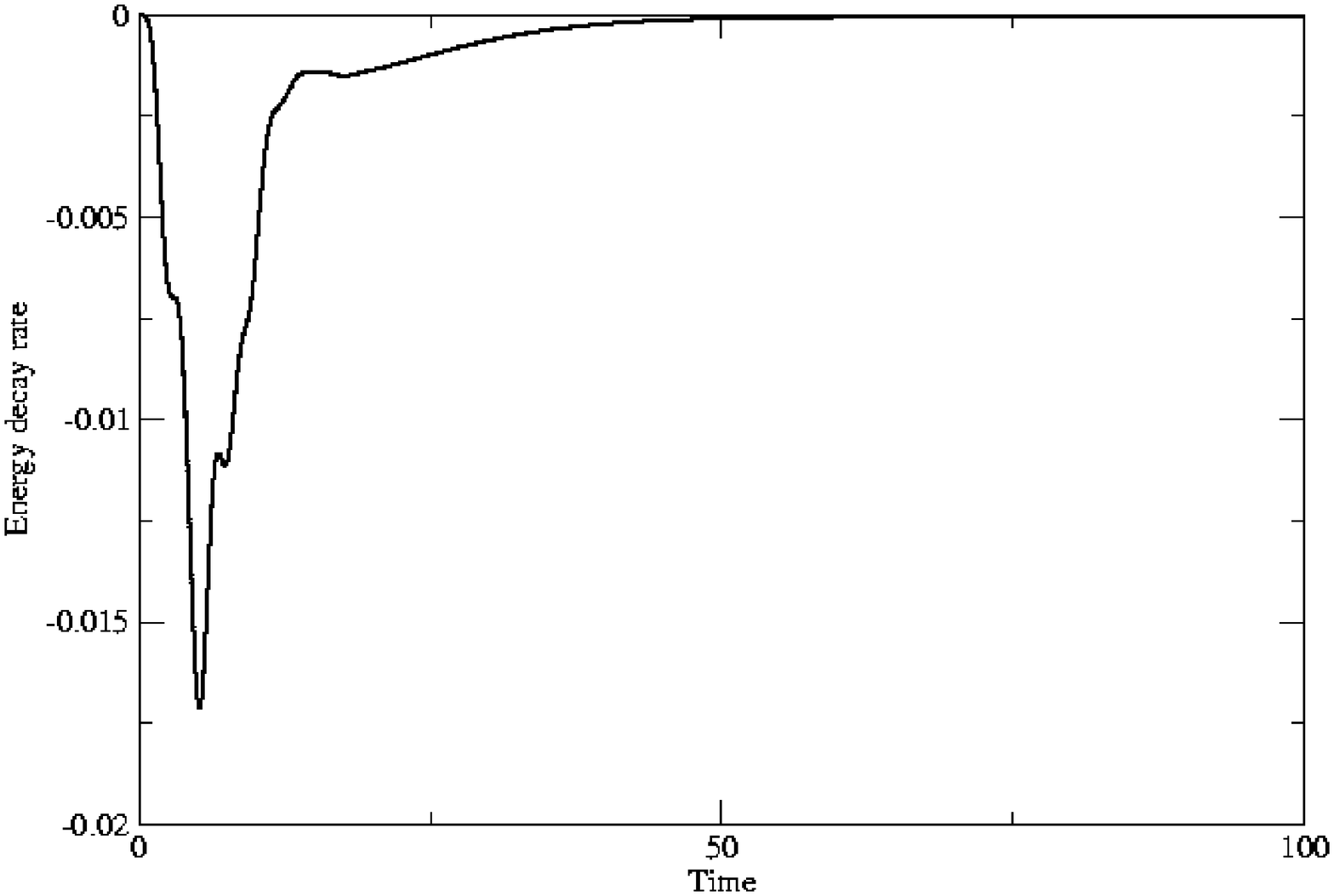,width=2.in}}
\caption{(a) Energy evolution for the first order model with $N=8^3$ modes. (b) Evolution of the energy decay rate. }
\label{fig_first}
\end{figure}

We continue our presentation of numerical results with the first order model. Figures \ref{fig_first}(a)  and (b) show the evolution of the energy and of the energy decay rate for a reduced system with $8^3$ modes.  The details of the numerical implementation are the same as in the case of the zeroth order model, including the truncation $t_0$ which is again set to $t_0=2.$ As discussed in Section \ref{zero}, not all expressions for the first order model are dealiased by construction. The necessary dealiasing is 
done through the 3/2 rule.

Inspection of the energy decay rate for the first and zeroth order model shows that the wavy structure 
of the energy decay has a shorter period for the first order model. This reflects in the energy decay rate 
as an increase in the number of spikes. Also, the slope of the energy decay (in log-log coordinates) has dropped to $\beta=-1.1379 \pm 0.0004,$ so the energy decays as $t^{\beta}.$ The theoretical estimate for isotropic decaying turbulence in the limit of infinite Reynolds number \cite{speziale} is $t^{-1}$ for 
the case of complete self-preservation. The notion of complete self-preservation means that the 
solution is self-similar across all scales, from zero to infinity. For solutions that are only partially self-preserved, the estimate is $t^{\gamma}$ with $\gamma < -1$. For periodic solutions in a box of finite size we cannot satisfy the complete self-preservation requirement, so that our estimate of a faster than $t^{-1}$ decay becomes more plausible. In addition, the Taylor-Green vortex problem is not isotropic, and 
we do not know yet how lack of isotropy can affect the energy decay.
During the simulations the magnitude of the zeroth term is one to two orders larger than the first order term.    

\begin{figure}
\centering
%\subfigure[]{\epsfig{file=fig_compare_energy.eps,width=2.in}}
\subfigure[]{\epsfig{file=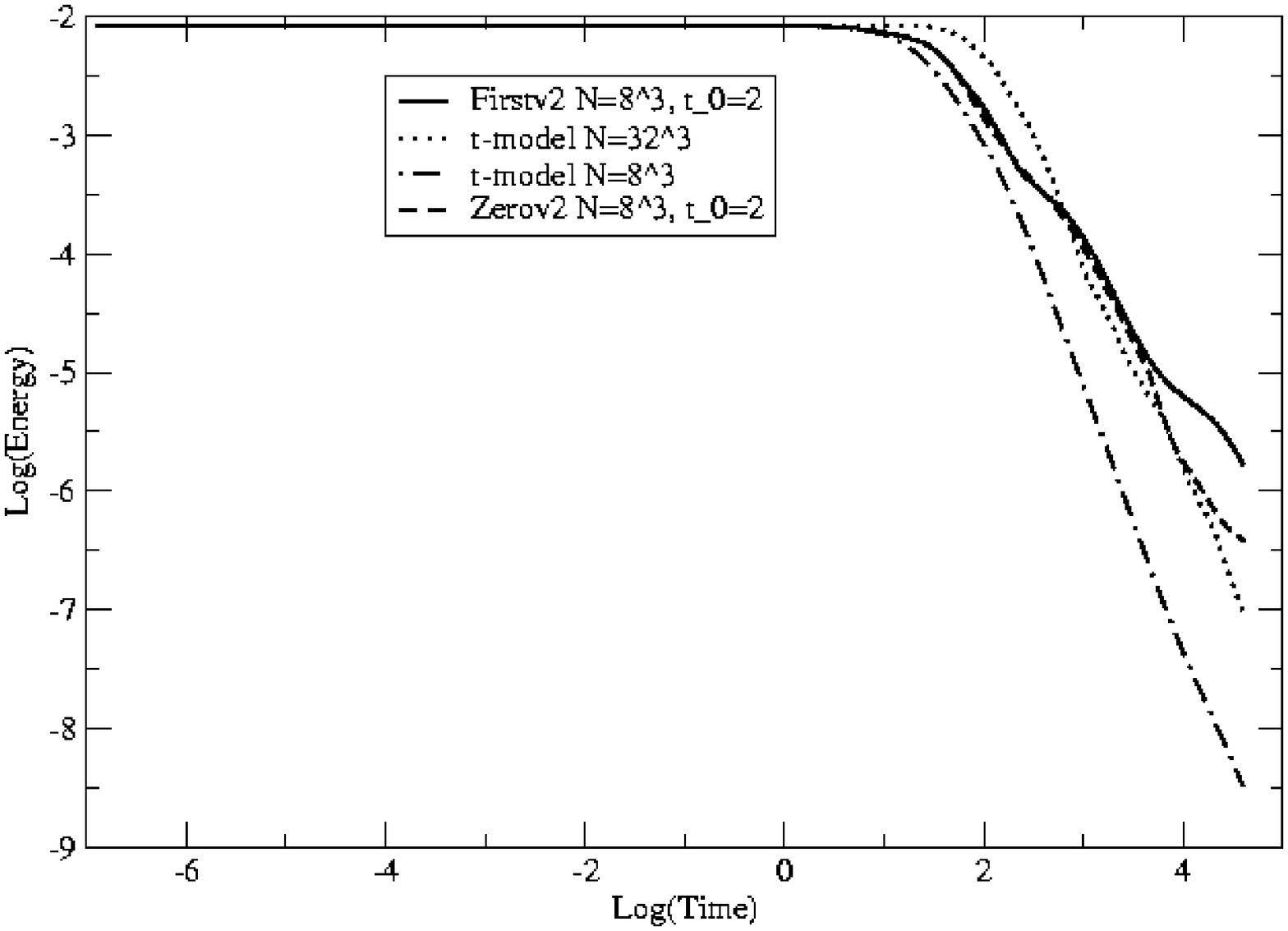,width=2.in}}
\qquad
%\subfigure[]{\epsfig{file=fig_compare_rate.eps,width=2.in}}
\subfigure[]{\epsfig{file=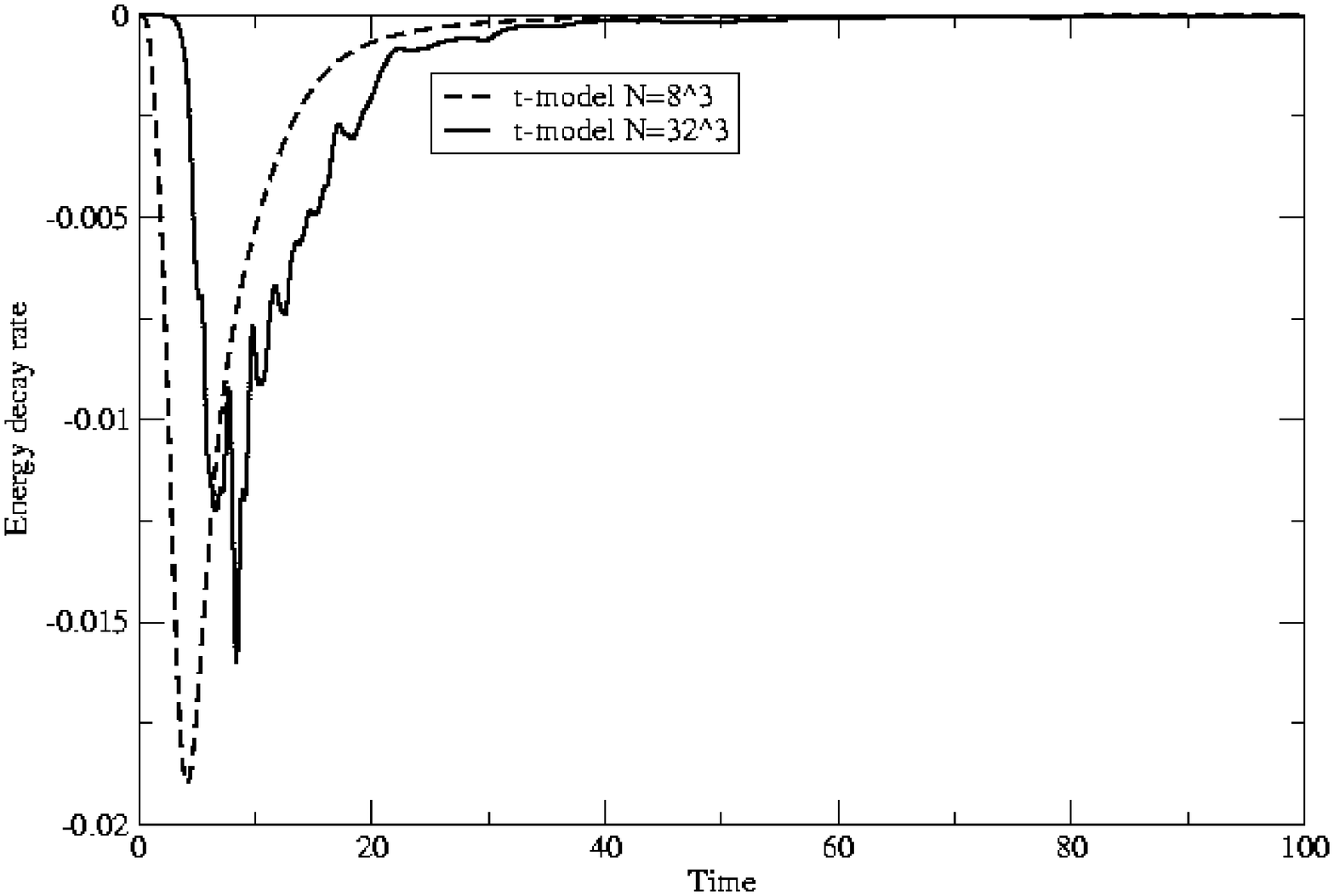,width=2.in}}
\caption{(a) Comparison of the energy evolution for the zeroth and first order models and the t-model with two different resolutions. (b) 
Rate of energy decay for t-model for reduced systems with $8^3$ and $32^3$ modes respectively.  }
\label{fig_compare}
\end{figure}

In Figure \ref{fig_compare}(a) we show how the energy decay of the zeroth and first order models with 
$8^3$ modes compare to the energy decay for the t-model with $8^3$ and $32^3$ modes. We see 
that as the resolution in the t-model is increased it resembles more the results of the models in this 
paper. The slopes of energy decay for the t-model are $-2.10$ and $-1.81$ for the cases with $8^3$ and $32^3$ modes respectively.
In Figure \ref{fig_compare}(b) we show the significant difference in the rate of energy decay 
for the t-model as we increase the number of resolved modes. The energy decay rate for $32^3$ modes 
exhibits organization of the energy decay in spikes, followed by periods of milder decay. It is 
important to note that the frequency of spikes in the $32^3$ case for the t-model is much higher than any other case presented here. We expect that the same increase in the frequency of spikes will appear also for the zeroth and first model presented here when we perform experiments with larger models. However, the point of this paragraph is that the use of more sophisticated models can reveal, at a lower resolution, effects that need a higher resolution if a lower order model is used. Also, it shows that the 
t-model albeit seemingly crude is on the right track if one can afford large calculations. Since the t-model is easier to implement than the models here, it appears as an interesting alternative to the more sophisticated models of this paper.

\subsection{Second order model}

\begin{figure}
\centering
%\subfigure[]{\epsfig{file=fig_second_energy.eps,width=2.in}}
\subfigure[]{\epsfig{file=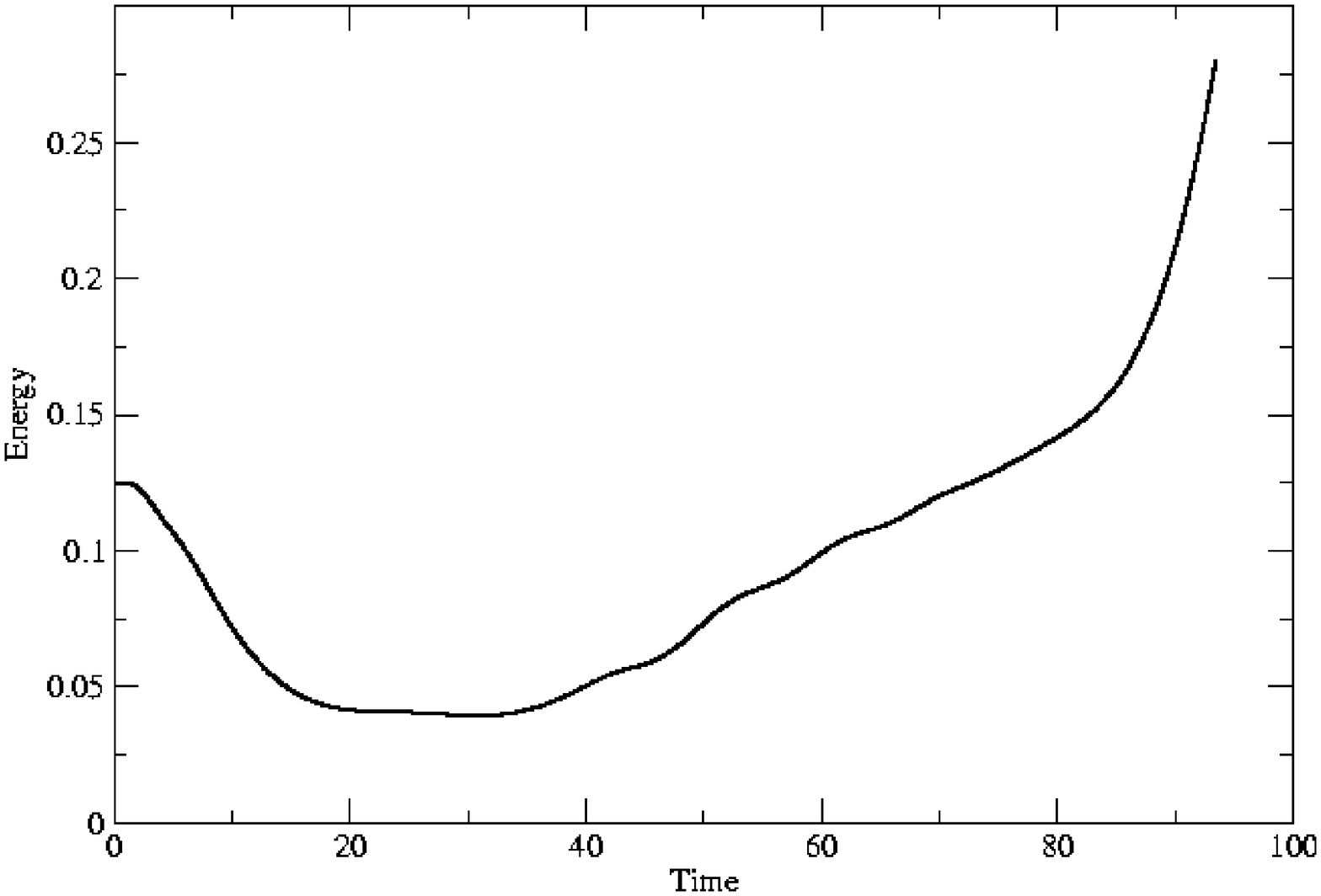,width=2.in}}
\qquad
%\subfigure[]{\epsfig{file=fig_second_rate.eps,width=2.in}}
\subfigure[]{\epsfig{file=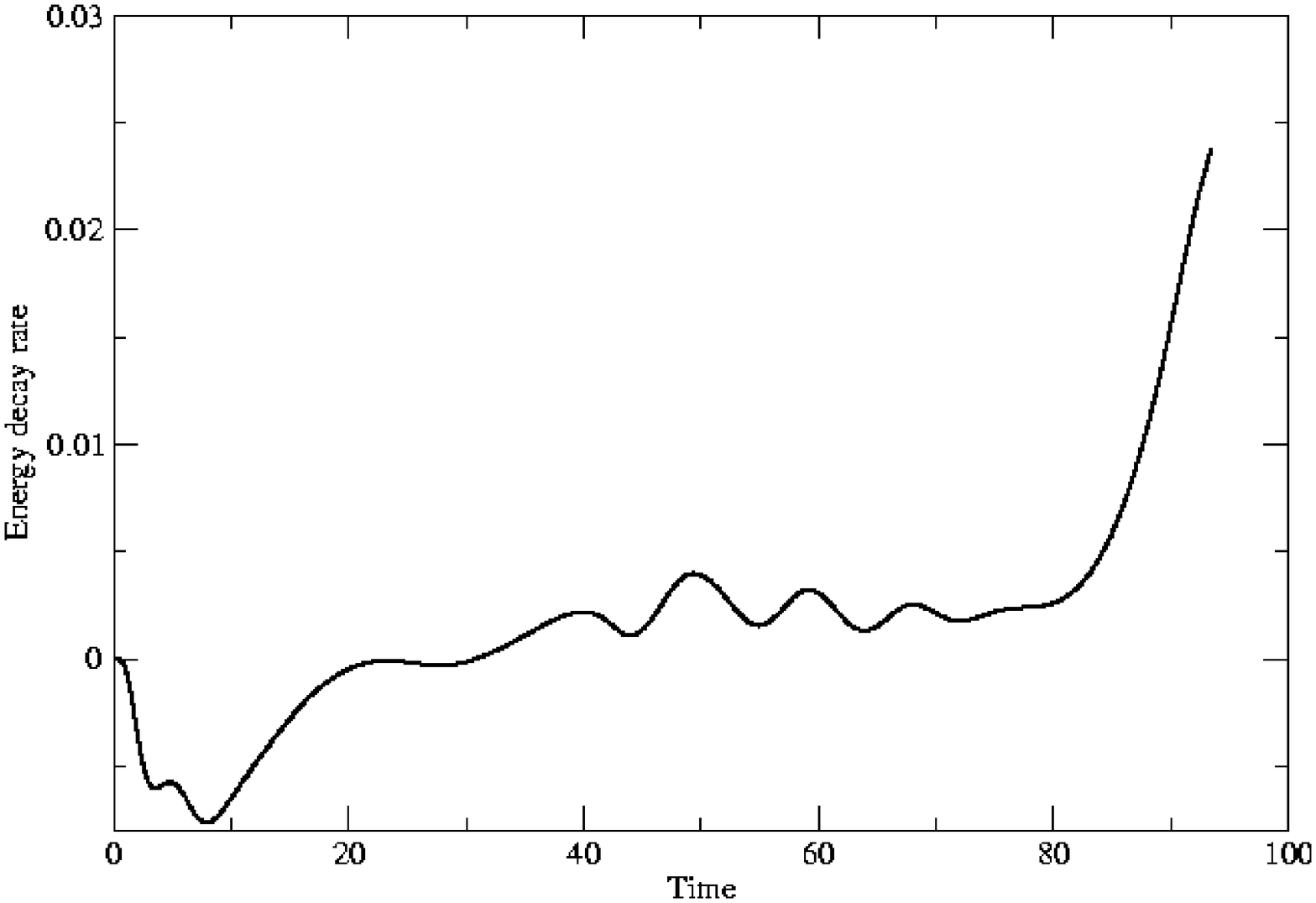,width=2.in}}
\caption{(a) Energy evolution for the second order model with $N=8^3$ modes. (b) Evolution of the energy decay rate. }
\label{fig_second}
\end{figure}

We conclude our presentation of the results for the reduced models with the second order model. Figures \ref{fig_second}(a) and (b) show the evolution of the energy and of the energy decay rate for a reduced system with $8^3$ modes.  As discussed in Section \ref{zero}, not all expressions for the second order model are dealiased by construction. The necessary dealiasing is done through the 3/2 rule.

The results shown are for $t_0=1.$ It is obvious that the reduced system is 
unstable. The energy cannot grow to a value larger than the initial one. No matter how small a value 
we tried for $t_0$ (from $t_0=2$ down to $t_0=0.01$), the instability is present. The instability should be related to 
a breakdown of the Taylor expansion. The energy decay rate is governed by the zeroth, first and second order terms. The magnitudes of the zeroth and second order terms are the main contributors, while the 
first order term (just as in the case of the first order model) is one order of magnitude smaller.

%%%%%%%%%%%%%%End of Section 5%%%%%%%%%%%%%%

\section{Conclusions}\label{conclusions}
We have presented a collection of reduced models for the 3D Euler equations based on the Taylor 
expansion of the orthogonal dynamics operator that appears in the Mori-Zwanzig formalism. The 
models up to second order were implemented and used to calculate the rate of energy decay for the 
Taylor-Green vortex problem. The results appear to be in good agreement with the theoretical 
estimates. The energy decay appears to be organized in "waves" of activity, i.e. 
alternating periods of fast and slow decay.

We presented a set of rules that can facilitate the recursive calculation of higher order models. The rules are based on the observation that the form of the terms appearing in the reduced models is determined by the form of the nonlinearity, the special properties of the projection operator used and the general properties of any projection operator.

The zero and first order models require a truncation of the range of integration for the 
memory term in order to produce solutions with decaying energy. However, there is no tuning needed. The results become better, the longer the range of integration, until the value of the range that leads to instability. Thus, trial and error is needed not to fit the results to some prescribed curve, but just to find when does the calculation become unstable.

On the other hand, the second 
order model does not produce solutions with decaying energy no matter how small is the range of integration for the memory term. The numerical results suggest that second and higher order terms in the Taylor expansion may be more profitably used in the construction of more elaborate approximations, e.g. Pad\'{e} approximants. Of 
course, before one starts thinking in terms of Pad\'{e} approximants, higher than second order terms in the Taylor expansion should be calculated and implemented. These higher order terms can potentially 
eliminate the instability that was found for the second order model. In addition, the numerical results 
suggest that a model including only a few low wavenumber modes should involve a long memory. 
This is in agreement with the main assumption behind the t-model presented in \cite{CHSS06}.

The problem of constructing reduced models for the Euler equations has been, and still is, a great 
challenge for scientific computing. The models proposed here should be considered a first step 
in deriving models directly from the equations without {\it ad hoc} approximations. They are based on numerical and physical observations about the behavior of the solution. The terms appearing in the reduced models can be efficiently implemented by the use of the FFT on appropriate arrays. This makes the incorporation of the models in existing pseudospectral algorithms rather straightforward. We plan to 
apply the models in a parallel setting which will allow a better assessment of the properties of 
the flow field that is predicted by the models.

\section{Acknowledgements} I am grateful to Profs G.I. Barenblatt, 
A.J. Chorin, O.H. Hald, Dr. Yelena Shvets and Mr. Jonathan Weare for many helpful discussions and comments.
This work was supported in part by the National Science 
Foundation under Grant DMS 04-32710, and by the Director,
Office of Science, Computational and Technology Research, U.S.\ Department of 
Energy under Contract No.\ DE-AC03-76SF000098.

\end{document}